\documentclass[12pt]{article}
\usepackage{amsmath, amsthm, amssymb,latexsym, enumerate, soul,ulem,cancel}
\usepackage{graphicx}
\usepackage{hyperref}
\usepackage{xcolor}
\usepackage{cite}
\usepackage{dsfont}
\usepackage{graphicx, caption, subcaption}

\newtheorem{theorem}{Theorem}[section]
\newtheorem{lemma}[theorem]{Lemma}
\newtheorem{corollary}[theorem]{Corollary}
\newtheorem{claim}[theorem]{Claim}

\newtheorem{conjecture}[theorem]{Conjecture}

\numberwithin{equation}{section}

\newcommand{\vp}{\varphi}
\newcommand{\D}{\Delta}

\DeclareMathOperator{\mad}{mad}

\title{Strong odd coloring of sparse graphs}

\author{Hyemin Kwon\thanks{
Korea Institute for Advanced Study,   85 Hoegi-ro, Dongdaemun-gu, Seoul 02455, Republic of Korea. 
\texttt{khmin1121@kias.re.kr}
} and
Boram Park\thanks{
Department of Mathematics Education, Seoul National University,  1 Gwanak-ro, Gwanak-gu, Seoul 08826,  Republic of Korea. 
\texttt{borampark@snu.ac.kr}
}
}
\date\today

\begin{document}

\maketitle
\begin{abstract}
An {\it odd coloring} of a graph $G$ is a proper coloring of $G$ such that for every non-isolated vertex $v$, there is a color appearing an odd number of times in $N_G(v)$. Odd coloring of graphs was studied intensively in recent few years.
In this paper, we introduce the notion of a strong odd coloring, as not only a  strengthened version of odd coloring, but also a relaxation of square coloring.  A {\it strong odd coloring} of a graph $G$ is a proper coloring of $G$ such that for every non-isolated vertex $v$, if a color appears in $N_G(v)$, then it appears an odd number of times in $N_G(v)$. We denote by
$\chi_{so}(G)$ the smallest integer $k$ such that $G$ admits a strong odd coloring with $k$ colors.
We prove that if $G$ is a graph with $mad(G)\le\frac{20}{7}$, then $\chi_{so}(G)\le \Delta(G)+4$, and the bound is tight.
We also prove that if $G$ is a $C_4$-free graph with $mad(G)\le\frac{30}{11}$, then $\chi_{so}(G)\le \Delta(G)+3$.
\end{abstract}

\section{Introduction}
All graphs in this paper are finite. Let $G$ be a graph.
For a vertex $v$, let $\deg_G(v)$ be the degree of $v$, $N_G(v)$ be the neighborhood of $v$, and $N_G[v]:=N_G(v)\cup \{v\}$.
Also, $\Delta(G)$ is the maximum degree of $G$.
The {\it girth} $g(G)$ of $G$ is the length of a shortest cycle in $G$ and the {\it maximum average degree} $\mad(G)$ of $G$ is the maximum of $\frac{2|E(H)|}{|V(H)|}$ over all non-empty subgraphs $H$ of $G$.

For a positive integer $k$, a \textit{proper $k$-coloring} of a graph is a function from the vertex set to $\{1,\cdots,k\}$ such that adjacent vertices receive different colors.
The minimum $k$ for which a graph $G$ has a proper $k$-coloring is the \textit{chromatic number} $\chi(G)$ of $G$.
A {\it square coloring} of a graph $G$ is a proper coloring of $G^2$, where $G^2$ is the graph obtained from $G$ by adding edges joining vertices at distance at most $2$ in $G$.
In 1977, Wegner conjectured the bound for $\chi(G^2)$ when $G$ is planar.

\begin{conjecture}[\cite{wegner1977}]\label{conj:wegner}
Let $G$ be a planar graph. Then
$$\chi(G^2)\leq
\begin{cases}
7 &\text{if $\Delta(G)=3$,}\\
\Delta(G)+5 &\text{if $4\leq \Delta(G) \leq 7$,}\\
\left\lfloor\frac{3\Delta(G)}{2}\right\rfloor+1 &\text{if $\Delta(G)\geq 8$.}
\end{cases}$$
\end{conjecture}

Thomassen\cite{thomassen2018} and Hartke et al.\cite{hartke2016} independently proved Conjecture~\ref{conj:wegner} when $\Delta(G)=3$, and it has attracted many researchers to study square coloring problems (see \cite{borodin2012,dvovrak2008list,la20212distance,cranston2022coloring,bu2012optimal,bu2018channel,deniz20232,dong20172,deniz2022improved,bu2016list} and Table~\ref{tab:square} for some results on planar graphs).

\begin{table}[h!]
    \centering
    \begin{tabular}{|c|c|c|c|c|c|c|c|}
         \hline
         &  $\D\ge 2$ & $\D\ge 4$ & $\D\ge 6$ & $\D\ge 9$ & $\D\ge 15$ & $\D\ge 22$ & $\D\ge 339$\\
         \hline
         $g\ge 5$ &  &  &  &  & $\D+5$\cite{bu2018channel} & $\D+4$\cite{deniz20232} & $\D+3$\cite{dong20172}\\
         \hline
         $g\ge 6$ & $\D+5$\cite{bu2012optimal} &  & $\D+4$\cite{deniz2022improved} & $\D+3$\cite{bu2016list} &  &  & \\
         \hline
         $g\ge 7$ &  &  & $\D+3$\cite{la20212distance} &  &  &  & \\
         \hline
         $g\ge 8$ &  & $\D+3$\cite{la20212distance} &  &  &  &  & \\
         \hline
    \end{tabular}
    \caption{Some known bounds on $\chi(G^2)$ for planar graphs $G$ with $\D(G)=\D$ and $g(G)=g$.}
    \label{tab:square}
\end{table}

The concept of odd coloring was introduced by Petru\v{s}evski and \v{S}krekovski~\cite{petrusevski2021colorings} as not only a strengthening of proper coloring and but also a weakening of square coloring.
For a positive integer $k$, an {\it odd $k$-coloring} of a graph $G$ is a proper $k$-coloring of $G$ such that every non-isolated vertex $v$ has a color appearing an odd number of times in $N_G(v)$.
The minimum $k$ for which $G$ has an odd $k$-coloring is the {\it odd chromatic number} of $G$, denoted by $\chi_o(G)$.
Ever since the first paper by Petru\v{s}evski and \v{S}krekovski~\cite{petrusevski2021colorings} appeared, there have been numerous papers~\cite{caro2022remarksodd,arXiv_ChChKwBo,cranston2022note,dujmovic2022odd,liu20221,petr2022odd,wang2022odd,cranston2022odd} studying various aspects of this new coloring concept across several graph classes.
Recently, Dai, Ouyang, and Pirot~\cite{dai2023new} introduced a strengthening of odd coloring.
For positive integers $h$ and $k$, an $h$-odd $k$-coloring is a proper $k$-coloring of a graph $G$ such that every vertex $v$ has $\min\{\deg_G(v),h\}$ colors each appearing an odd number of times in $N_G(v)$.
The $\chi_o^h(G)$ is the minimum $k$ for which $G$ has an $h$-odd $k$-coloring.
It is clear that $\chi_o(G)\le \chi_o^h(G) \le \chi(G^2)$ by definitions.

\begin{figure}[h!]
    \centering
    \includegraphics[width=0.3\textwidth,page=1]{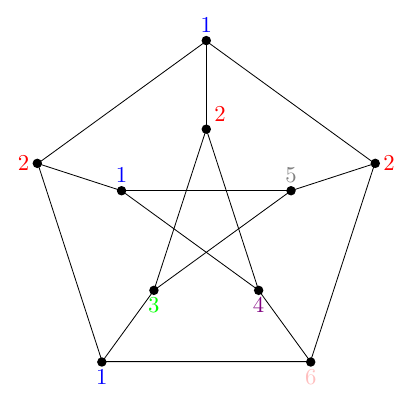}
    \caption{An SO $6$-coloring of the Petersen graph.}
    \label{fig:petersen}
\end{figure}

In this paper, we introduce a new concept that is a new generalization of odd coloring.
For a positive integer $k$, a \textit{strong odd $k$-coloring} (an SO $k$-coloring for short) of a graph $G$ is a proper $k$-coloring of $G$ such that for every non-isolated vertex $v$, if a color appears in $N_G(v)$, then it appears an odd number of times in $N_G(v)$.
The \textit{strong odd chromatic number} of a graph $G$, denoted by $\chi_{so}(G)$, is the minimum $k$ such that $G$ has an SO $k$-coloring.
See Figure~\ref{fig:petersen} for an example.
Like as odd coloring, strong odd coloring does not have hereditary, that is, there are graphs $H$ and $G$ such that $H\subset G$ and $\chi_{so}(H)>\chi_{so}(G)$.
For example, $\chi_{so}(K_{1,2})=3$ and $\chi_{so}(K_{1,3})=2$ while $K_{1,2}\subset K_{1,3}$.

Every graph has a strong odd coloring since a square coloring is also a strong odd coloring.
Moreover, it is clear from the definitions that for every graph $G$,
\[ \chi_o(G)\le \chi_{so}(G)\le \chi(G^2)\le (\D(G))^2+1.\]
When $G$ is claw-free, that is, $G$ has no $K_{1,3}$ as an induced subgraph, it holds that $\chi_{so}(G)=\chi(G^2)$.
This can be seen because when $G$ is claw-free, for a vertex $v$ of $G$, no matter how three neighbors of $v$ are chosen, two of them are adjacent.
Therefore, the three neighbors must all have different colors by the definition of a strong odd coloring.
However, $\chi(G^2)-\chi_{so}(G)$ could be arbitrary large. For example, if $G=K_{m,n}$, then $\chi(G^2)=m+n$ and $\chi_{so}(G)\le 4$.
We remark that whereas $\chi_o(G)\le 2\Delta(G)+1$ by \cite{caro2022remarksodd}, a linear bound of $\chi_{so}(G)$ in terms of $\Delta(G)$ cannot be obtained in general.
For example, if $G=K_n \square K_n$, then $\Delta(G)=2n-2$, the diameter of $G$ is two, and $G$ is claw-free, and thus $\chi_{so}(G)=\chi(G^2)=|V(G)|=n^2=\frac{(\D(G)+2)^2}{4}$.

\begin{figure}[h!]
    \centering
    \includegraphics[width=0.25\textwidth,page=2]{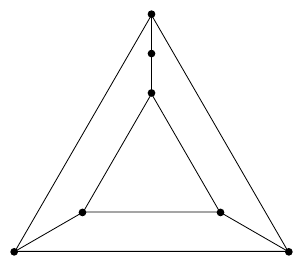}
    \caption{A subcubic planar graph $G$ with $\chi(G^2)=\chi_{so}(G)=7$}
    \label{fig:planarcubic7}
\end{figure}

Regarding to subcubic planar graphs $G$, we have $\chi_{so}(G)\le \chi(G^2)\le 7$, the bound is tight since there is a subcubic planar graph $G$ such that $\chi_{so}(G)=\chi(G^2)=7$. 
Considering the graph in Figure~\ref{fig:planarcubic7}, it has strong odd chromatic number $7$ and it has maximum average degree $\frac{20}{7}$. This provides a tight example to the following main result.

\begin{theorem}\label{thm:soD+4}
If $G$ is a graph with $mad(G)\le\frac{20}{7}$, then $\chi_{so}(G)\le \Delta(G)+4$.
\end{theorem}

In \cite{la20212distance}, it was shown that $\chi(G^2)\le \Delta(G)+3$ if $G$ is a graph with $\D(G)\ge 4$ and $\mad(G)< \frac{8}{3}$ and therefore, we have the same bound for $\chi_{so}(G)$. 
We find a condition of a graph $G$ satisfying $\chi_{so}(G)\le \Delta(G)+3$.
\begin{theorem}\label{thm:so6coloring}
If $G$ is a $C_4$-free graph with $mad(G)\le\frac{30}{11}$, then $\chi_{so}(G)\le \Delta(G)+3$.
\end{theorem}


Since $(\mad(G)-2)(g(G)-2)<4$ for every planar graph $G$,  we obtain the following corollary from Theorems~\ref{thm:soD+4}, and \ref{thm:so6coloring}.  See also Table \ref{tab:square} to compare the results on square coloring of planar graphs with large maximum degree and large girth.
\begin{corollary}\label{cor:soD+4}
Let $G$ be a planar graph.
\begin{itemize}
    \item[\rm(i)] If $g(G)\ge 7$, then $\chi_{so}(G)\le \Delta(G)+4$.
    \item[\rm(ii)] If $g(G)\ge 8$, then $\chi_{so}(G)\le \Delta(G)+3$.
\end{itemize}
\end{corollary}

The paper is organized as follows. In Section 2, we collect some basic observations on $\chi_{so}(G)$.
Section 3 proves Theorem~\ref{thm:soD+4}, and Section 4 proves Theorem~\ref{thm:so6coloring}.
We make some final remarks in Section 5.

\section{Preliminaries}

For $S\subset V(G)$, let $G-S$ denote the graph obtained from $G$ by deleting the vertices in $S$.
If $S=\{x\}$, then denote $G-S$ by $G-x$.
For a positive integer $k$, we use $k$-vertex (resp. $k^+$-vertex and $k^-$-vertex) to denote a vertex of degree $k$ (resp. at least $k$ and at most $k$).
We also use $k$-neighbor (resp. $k^+$-neighbor and $k^-$-neighbor) of a vertex $v$ to denote a $k$-vertex (resp. $k^+$-vertex and $k^-$-vertex) that is a neighbor of $v$.
For positive integers $k$ and $d$, we call a $k$-vertex with at least $d$ $2$-neighbors  a \textit{$k_d$-vertex}.
We also use \textit{$k_d$-neighbor} of a vertex $v$ to call a neighbor of $v$ that is  a $k_d$-vertex.
 
For an integer $d\ge 2$, let $A_1,\ldots,A_{d}$ be nonempty sets. A sequence $(c_1,\ldots, c_d)$ is called a \textit{system of odd representative} for $(A_1,\ldots,A_{d})$ if it satisfies the following:
\begin{itemize}
\item[(i)] $c_i\in A_i$ for each $1\le i\le d$, and
\item[(ii)] if $c=c_t$ for some $t \in \{1,\ldots,d\}$, then $c$ appears an odd number of times in $(c_1,\ldots, c_d)$.
\end{itemize}

\begin{lemma}\label{lem:even}
For an integer $d\ge 2$, let $A_1,\ldots,A_{d-1}$ be nonempty sets of size at least $2$.
Then for every $\alpha\in A_d$, there is a system of odd representative $(\beta_1,\ldots, \beta_{d-1},\alpha)$ for $(A_1,\ldots,A_{d})$.
\end{lemma}

\begin{proof}
We will prove the lemma using induction on $d$.
It is enough to assume that $|A_i|= 2$ for each $i\in\{1,\ldots,d-1\}$. Fix $\alpha\in A_d$.
It is clear that if $d=2$, we can choose $\beta_1\in A_1\setminus\{\alpha\}$ to find a system of odd representative $(\beta_1, \alpha)$ for $(A_1,A_2)$.
Let $d=3$.
If both $A_1$ and $A_2$ have $\alpha$ as an element, then we choose $\beta_1=\beta_2=\alpha$ and it is a system of odd representative $(\beta_1, \beta_2,\alpha)$ for $(A_1,A_2,A_3)$.
So we assume $\alpha\not\in A_1$, without loss of generality.
Then we choose $\beta_2\in A_2\setminus \{\alpha\}$, and $\beta_1\in A_1\setminus \{\beta_2\}$.
Thus, we find a system of odd representative $(\beta_1, \beta_2,\alpha)$ for $(A_1,A_2,A_3)$.

Let $d\ge 4$.
For given $A_1,\ldots, A_{d-1}$, suppose there are two sets containing $\alpha$.
Without loss of generality, we call the sets $A_{d-2}$ and $A_{d-1}$.
Then by inductive hypothesis, there is a system of odd representative $(\beta_1,\ldots, \beta_{d-3},\alpha)$ for $(A_1,\ldots, A_{d-3},A_d)$.
By choosing $\beta_{d-2}=\beta_{d-1}=\alpha$, we can find a system of odd representative $(\beta_1,\ldots, \beta_{d-1},\alpha)$ for $(A_1,\ldots, A_{d-1},A_d)$.
Now we consider the case when there is at most one set containing $\alpha$ for given $A_1,\ldots, A_{d-1}$.
Without loss of generality, let $\alpha\not\in A_1\cup\cdots \cup A_{d-2}$.
Then, we choose $\beta_{d-1}\in A_{d-1}\setminus\{\alpha\}$.
By inductive hypothesis, there is a system of odd representative $(\beta_1,\ldots, \beta_{d-2},\beta_{d-1})$ for $(A_1,\ldots, A_{d-2},A_{d-1})$.
Note that $\beta_i\neq \alpha$ for all $i$, so there is a system of odd representative $(\beta_1,\ldots, \beta_{d-1},\alpha)$ for $(A_1,\ldots, A_{d-1},A_{d})$.
\end{proof}

Throughout this paper, we use $\mathcal{C}$ to denote the set of all $k$ colors, when we consider {a} $k$-coloring of a graph $G$.
For a (partial) coloring $\vp$ of a graph $G$, let $\vp(X)=\{\vp(v)\mid v\in X\}$ for $X\subset V(G)$.

In each proof of following lemmas,  we will start a proof with defining a nonempty subset $S\subset V(G)$.
We always let $G'=G-S$, $G''=G^2-S$, and $\vp$ be  an SO $(\D(G)+c)$-coloring of $G'$, if it exists, for a positive integer $c$.
We also define \[A_{\vp}(v)=\mathcal{C}\setminus \vp(N_{G''}(v)).\]
We  call an element of $A_{\vp}(v)$ an \textit{available color} of $v$ under $\vp$.
Note that
\begin{eqnarray}\label{eq:Av}
|A_{\vp}(v)|&=& |\mathcal{C}|-|\vp(N_{G''}(v))|\ge \D(G)+c-\deg_{G''}(v).
\end{eqnarray}
When we color the vertices of $S\setminus S^*$ for some $S^*\subset S$, we often use the same symbol $\vp$ to denote the resulting coloring, that is a coloring of $G-S^*$.

\begin{lemma}\label{lem:G-x:2vtex}
For an integer $c\ge 3$, let $G$ be a minimal graph with respect to $|V(G)|$ such that $G$ has no SO $(\D+c)$-coloring, where $\Delta:=\Delta(G)$.
For a $2$-vertex $x$, there is no SO $(\D+c)$-coloring $\vp$ of $G-x$ such that there is an available color of $x$ under $\vp$ and the neighbors of $x$ get distinct colors.
\end{lemma}
\begin{proof} Suppose that such coloring $\vp$ exists. Then  by coloring $x$ with a color in $A_{\vp}(x)$, we obtain an SO $(\D+c)$-coloring of $G$. It is a contradiction to the fact that $G$ has no SO $(\D+c)$-coloring.
\end{proof}

 In all figures, we represent a vertex that has all incident edges in the figure as a filled vertex, and a hollow vertex may not have all incident edges in the figure.

\begin{lemma}\label{lem:D+c}
For an integer $c\ge 3$, let $G$ be a minimal graph with respect to $|V(G)|$ such that $G$ has no SO $(\D+c)$-coloring, where $\Delta:=\Delta(G)$.
Then the following do not appear in $G$:
\begin{enumerate}[\rm(i)]
    \item a $1^-$-vertex,
    \item a triangle $v_1v_2v_3v_1$ such that $v_2$ is a $(c+1)^-$-vertex and $v_3$ is a $2$-vertex,
    \item two vertices with three common $2$-neighbors,
    \item 
 a $d_{d-1}$-vertex with only $(\D+c-d)^-$-neighbors for $2\le d\le \Delta$,  when either $c\ge 4$ or $c=3$ and $G$ is $C_4$-free.
\end{enumerate}
\end{lemma}

\begin{figure}[h!]
    \centering\includegraphics[width=0.6\columnwidth,page=3]{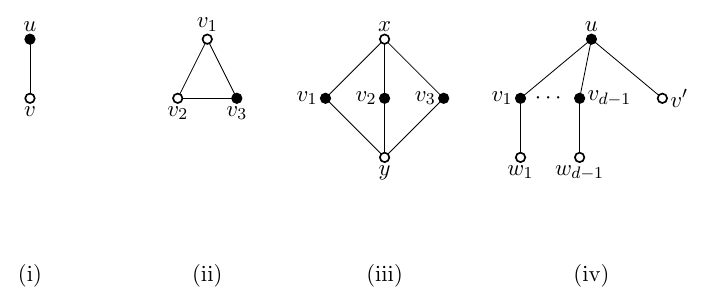}
    \caption{Illustrations of Lemma~\ref{lem:D+c}.}
    \label{fig:lem:D+c}
\end{figure}

\begin{proof}
In each proof, we consider an SO $(\Delta+c)$-coloring $\vp$ of $G'=G-S$ for some $\emptyset\neq S\subset V(G)$, and then we finish the proof by coming up with an SO $(\Delta+c)$-coloring of $G$, which is a contradiction.
See Figure~\ref{fig:lem:D+c}.

(i) It is clear that $G$ has no isolated vertex. Suppose to the contrary that $G$ has a $1$-vertex $u$ with neighbor $v$.
Let $S=\{u\}$.
Since $|A_{\varphi}(u)|\ge c$ by \eqref{eq:Av}, we color $u$ with a color in $A_{\varphi}(u)$, which gives an SO $(\Delta+c)$-coloring of $G$.

(ii) Suppose to the contrary that $G$ has a triangle $v_1v_2v_3v_1$ such that $v_2$ is a $(c+1)^-$-vertex and $v_3$ is a $2$-vertex.
For an SO $(\D+c)$-coloring $\vp$ of $G-v_3$, $|A_{\vp}(v_3)|\ge \D+c-(\D+c-1)\ge 1$ by \eqref{eq:Av}. Since $v_1v_2\in E(G')$, $\vp(v_1)\neq \vp(v_2)$.  It is a contradiction by Lemma~\ref{lem:G-x:2vtex}.

(iii) Suppose to the contrary that $G$ has two vertices $x$ and $y$ with three common $2$-neighbors $v_1$, $v_2$, and $v_3$.
Let $S=\{v_2,v_3\}$.
Then color $v_2$ and $v_3$ with the same color of $v_1$ to obtain an SO $(\Delta+c)$-coloring of $G$.

(iv) Suppose to the contrary that $G$ has a $d_{d-1}$-vertex $u$ with only $(\D+c-d)^-$-neighbors for some $2\le d\le \D$.
Let $v_1,\ldots, v_{d-1}$ be $2$-neighbors of $u$ and $v'$ be the other neighbor of $u$.
For each $i\in\{1,\ldots,d-1\}$, let $w_i$ be the neighbor of $v_i$ other than $u$.
Note that $v_i\neq w_j$ for each $i,j\in\{1,\ldots,d-1\}$ by (ii).
Let $S=\{u,v_1,\ldots,v_{d-1}\}$.
Note that $|A_{\vp}(u)|\ge \Delta+c -(d-1+\D+c-d)\ge 1$, and  $|A_{\vp}(v_i)|\ge \D+c-(\D+1)=c-1\ge 2$  for each $i\in\{1,\ldots,d-1\}$ by \eqref{eq:Av}.
We color $u$ with a color $\gamma$ in $A_{\vp}(u)$. Let $X$ be the set of all pairs $(i,j)$ such that $w_i=w_j$ and $i<j$. Note that by (iii), each $i\in\{1,\ldots, d-1\}$ appears at most once as an entry of an element of $X$. Moreover, when $(i,j)\in X$, it holds that $A_{\vp}(v_i)=A_{\vp}(v_j)$ and $|A_\vp(v_i)|=|A_\vp(v_j)|\ge c$ by \eqref{eq:Av}.

Suppose that $c=3$ and 
$G$ is $C_4$-free. Then $X= \emptyset$.
For each $i\in\{1,\ldots,d-1\}$, let 
\[ A_i:=\begin{cases}
    (A_{\vp}(v_i)\cup\{\vp(v')\})\setminus\{\gamma\}  &  \text{if }|A_{\vp}(v_i)|=2, \\
A_{\vp}(v_i)\setminus\{\gamma\} &\text{ otherwise}.
\end{cases}\]
Let $A_d:=\{\varphi(v')\}$.
Since $v'$ is a neighbor of $v_i$ in $G^2$, $\vp(v')\not\in A_{\vp}(v_i)$.
It follows that $|A_i|\ge 2$ for each $i\in\{1,\ldots,d-1\}$.
By applying Lemma~\ref{lem:even}, there is a system of odd representative $(\beta_1,\ldots,\beta_{d-1},\alpha)$ for $(A_1,\ldots,A_{d-1},A_d)$, where $\alpha=\varphi(v')$.
Then, coloring $v_i$ with $\beta_i$ results in an SO $(\Delta+c)$-coloring of $G$.

Suppose $c\ge 4$.
For each $i\in\{1,\ldots,d-1\}$, let
\[ A_i:=\begin{cases}
    (A_{\vp}(v_i)\cup\{\vp(v')\})\setminus\{\gamma\}  &  \text{if }|A_{\vp}(v_i)|=4 \text{ and }(i,j)\in X, \\
A_{\vp}(v_i)\setminus\{\gamma\} &\text{otherwise}.
\end{cases}\]
Let $A_d:=\{\varphi(v')\}$.
Since $v'$ is a neighbor of $v_i$ in $G^2$, $\vp(v')\not\in A_{\vp}(v_i)$.
It follows that $|A_i| \ge 2$ for each $i\in\{1,\ldots,d-1\}$. 
If $(i,j)\in X$, then $X\neq \emptyset$ and so $|A_i|\ge 4$, and then we redefine $A_i:=A'_i$ and $A_j:=A'_j$ for some disjoint subsets $A'_i,A'_j$ of $A_i$ of size two.
By applying Lemma~\ref{lem:even}, there is a system of odd representative $(\beta_1,\ldots,\beta_{d-1},\alpha)$ for $(A_1,\ldots,A_{d-1},A_d)$, where $\alpha=\varphi(v')$.
Then, 
coloring $v_i$ with $\beta_i$ results in an SO $(\Delta+c)$-coloring of $G$. Note that, when $w_i=w_j$, the definition of $A_i$ and $A_j$ guarantees that every color in the neighbors of $w_i$ appears odd number of times. 
\end{proof}

We finish the section by stating one famous theorem in graph coloring, called Brooks' Theorem.
\begin{theorem}[\cite{1941Brooks}]\label{thm:brooks}
For a graph $G$, $\chi(G)\le \D(G)+1$.
If $G$ has no component that is a complete graph or an odd cycle, then $\chi(G)\le \D(G)$.
\end{theorem}

\section{Strong odd $(\D+4)$-coloring}
In this section, we prove Theorem~\ref{thm:soD+4}.
Let $G$ be a minimal counterexample to Theorem~\ref{thm:soD+4}. In each proof in the following, we always start with defining a nonempty set $S\subset V(G)$,  and $\vp$ is an SO $(\Delta+{4})$-coloring  of $G'=G-S$. We end up with an SO $(\Delta+{4})$-coloring of $G$.

The following is a list of reducible configurations, structures that never appear in $G$.
The configurations \textbf{[C\ref{rc:D+4:1}]}-\textbf{[C\ref{rc:D+4:5}]} are utilized in the final step to reach a contradiction.
\begin{enumerate}[\bfseries {[}C1{]}]
    \item\label{rc:D+4:1} A $1^-$-vertex. \vspace{-0.25cm}
    \item\label{rc:D+4:2} A $d_{d-1}$-vertex for all $2\le d\le 4$. \vspace{-0.25cm}
    \item\label{rc:D+4:add} A $d_d$-vertex for all $d\ge 5$. \vspace{-0.25cm}
    \item\label{rc:D+4:3} Two adjacent $3_1$-vertices. \vspace{-0.25cm}
    \item\label{rc:D+4:4} A $3$-vertex with two $3_1$-neighbors. \vspace{-0.25cm}
    \item\label{rc:D+4:5} A $4_2$-vertex with a $3_1$-neighbor.
\end{enumerate}
Note that \textbf{[C\ref{rc:D+4:1}]}-\textbf{[C\ref{rc:D+4:add}]} do not appear in $G$ by Lemma~\ref{lem:D+c}~(i) and (iv).

Suppose that a $2$-vertex $x$ has a $3$-neighbor, and  $\vp$ is an SO $(\D+4)$-coloring of $G-x$. Since $|\mathcal{C}|=\D+4$, it holds that $A_{\vp}(x)\neq\emptyset$ by \eqref{eq:Av}. Thus by Lemma~\ref{lem:G-x:2vtex}, the colors of the neighbors of   $x$ cannot be distinct colors. We often omit this explanation when we apply Lemma~\ref{lem:G-x:2vtex}.

\begin{lemma}\label{lem:D+4}
The graph $G$ has no two adjacent $3_1$-vertices. {\rm\textbf{[C\ref{rc:D+4:3}]}}
\end{lemma}
\begin{figure}[h!]
    \centering
    \includegraphics[width=0.2\columnwidth,page=4]{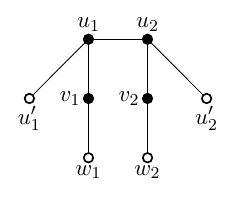}
    \caption{An illustration of Lemma~\ref{lem:D+4}.}
    \label{fig:lem:D+4}
\end{figure}
\begin{proof}
Suppose to the contrary that $G$ has two adjacent $3_1$-vertices $u_1$ and $u_2$.
For each $i\in\{1,2\}$, let $v_i$ be the $2$-neighbor of $u_i$, $u'_i$ be the neighbor of $u_i$ other than $v_i$ and $u_{3-i}$, and $w_i$ be the neighbor of $v_i$ other than $u_i$.
See Figure~\ref{fig:lem:D+4}.
Note that $v_1\neq v_2$ by Lemma~\ref{lem:D+c}~(ii). Thus $u_1$, $u_2$, $v_1$, and $v_2$ are distinct, and we let $S=\{u_1,u_2,v_1,v_2\}$.
In addition, $w_i\neq u_{3-i}$ (equivalently, $u'_{3-i}\neq v_i$) and $w_i\neq v_{3-i}$ for each $i\in\{1,2\}$ by {\rm\textbf{[C\ref{rc:D+4:2}]}}. 
Thus  $S\cap \{u'_1,u'_2,w_1,w_2\}=\emptyset$.

Note that $|A_{\vp}(u_i)|\ge 2$ and $|A_{\vp}(v_i)|\ge 3$ for each $i\in\{1,2\}$ by \eqref{eq:Av}.
We color $u_1$, $u_2$, $v_1$ with a color in $A_{\vp}(u_1)$,  $A_{\vp}(u_2)$, $A_{\vp}(v_1)$, respectively, so that colors of $u_1$, $u_2$, $v_1$ are distinct. Then the resulting coloring is an SO $(\D+4)$-coloring of $G-v_2$, which contradicts to Lemma~\ref{lem:G-x:2vtex}.
\end{proof}

\begin{lemma}\label{lem:vertexdistinct}
The graph $G$ has no path $v_1v_2v_3v_4$ such that $v_1$ is a $2$-vertex, $v_i$ is a $3$-vertex for each $i\in\{2,3,4\}$, where $v_4$ is a common neighbor of $v_3$ and $v_j$ for some $j\in \{1,2\}$.
\end{lemma}
\begin{figure}[h!]
    \centering
    \includegraphics[width=0.4\columnwidth,page=5]{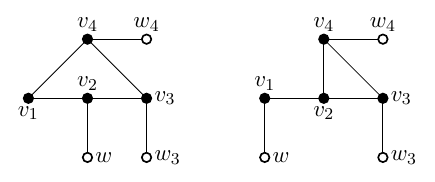}
    \caption{Illustrations of Lemma~\ref{lem:vertexdistinct}.}
    \label{fig:lem:vertexdistinct}
\end{figure}
\begin{proof}
Suppose to the contrary that $G$ has a path $v_1v_2v_3v_4$ such that $v_1$ is a $2$-vertex and $v_i$ is a $3$-vertex for each $i\in\{2,3,4\}$, where $v_4$ is a common neighbor of $v_3$ and either $v_1$ or $v_2$.
Other neighbors of $v_i$'s are labeled as Figure~\ref{fig:lem:vertexdistinct}.
Let $S=\{v_1,v_2,v_3,v_4\}$.
By Lemma~\ref{lem:D+c}~(ii), $S\cap\{w,w_3,w_4\}=\emptyset$.
Note that $|A_{\vp}(v_1)|\ge 4$, $|A_{\vp}(v_2)|\ge 3$, $|A_{\vp}(v_3)|\ge 2$, and $|A_{\vp}(v_4)|\ge 3$ by \eqref{eq:Av}.
Then it is possible to color $v_i$ with a color in $A_{\vp}(v_i)$ for each $i\in\{1,2,3,4\}$ so that the colors of $v_i$'s are distinct.
It results in an SO $(\D+4)$-coloring of $G$.
\end{proof}

\begin{lemma}\label{lem:D+4:two3_1}
The graph $G$ has no $3$-vertex with two $3_1$-neighbors. {\rm\textbf{[C\ref{rc:D+4:4}]}}
\end{lemma}
\begin{figure}[h!]
    \centering
    \includegraphics[width=0.25\columnwidth,page=6]{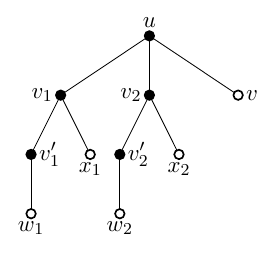}
    \caption{An illustration of Lemma~\ref{lem:D+4:two3_1}.}
    \label{fig:lem:D+4:two3_1}
\end{figure}
\begin{proof}
Suppose to the contrary that $G$ has a $3$-vertex $u$ with two $3_1$-neighbors $v_1$ and $v_2$.
We follow the labeling of the vertices as Figure~\ref{fig:lem:D+4:two3_1}.
Let $S=\{u,v_1,v_2,v'_1,v'_2\}$. 
By Lemma~\ref{lem:vertexdistinct}, $v'_1\neq v'_2$, and so $S$ has five distinct vertices.
Note that for each $i\in\{1,2\}$, $v_i\notin \{w_{3-i},x_{3-i}\}$ by Lemma~\ref{lem:vertexdistinct}.
Also, $u\neq w_i$, $v'_i\neq v$ by Lemma~\ref{lem:D+c}~(ii), and $v'_i\notin \{w_{3-i},x_{3-i}\}$ by Lemma~\ref{lem:D+c}~(iv).
Thus $S\cap\{ w_1,w_2,x_1,x_2,v\}=\emptyset$.
Note that $|A_{\vp}(u)|\ge 2$, $|A_{\vp}(v_i)|\ge 2$, and $|A_{\vp}(v'_i)|\ge 3$ for each $i\in\{1,2\}$ by \eqref{eq:Av}.

First, we color $u$ with a color $\gamma$ in $A_{\vp}(u)$.
For simplicity, let $A_i:=(A_{\vp}(v_i)\cup\{\vp(v)\})\setminus\{\gamma\}$ if $|A_{\vp}(v_i)|=2$, and $A_i:=A_{\vp}(v_i)\setminus\{\gamma\}$ otherwise.
Since $\vp(v)\notin A_{\vp}(v_i)$, $|A_i| \ge 2$ for each $i\in\{1,2\}$.
Let $A_3=\{\vp(v)\}$.
If $x_1\neq x_2$, then we apply Lemma~\ref{lem:even} to obtain a system of odd representative $(\beta_1,\beta_2,\alpha)$ for $(A_1,A_2,A_3)$, where $\alpha=\vp(v)$.
If $x_1=x_2$, then $|A_{\vp}(v_i)|\ge 3$ for each $i\in\{1,2\}$ by \eqref{eq:Av}, and so we choose $\beta_i$ in $A_{\vp}(v_i)\setminus\{\gamma\}$ for each $i\in\{1,2\}$ so that $\beta_1\neq \beta_2$.
Then we color $v_i$ with $\beta_i$ in $A_i$ for each $i\in\{1,2\}$.
We denote this coloring of $G-\{v'_1,v'_2\}$ by $\vp$ again. 
Then $|A_{\vp}(v'_i)|\ge 1$ for each $i\in\{1,2\}$  by \eqref{eq:Av} and so we color $v'_1$ with a color in $A_{\vp}(v'_1)$.
Then the resulting coloring is an SO $(\D+4)$-coloring of $G-v'_2$, which contradicts to Lemma~\ref{lem:G-x:2vtex}.
\end{proof}

\begin{lemma}\label{lem:D+4:4two2}
The graph $G$ has no $4_{2}$-vertex with a $3_1$-neighbor. {\rm\textbf{[C\ref{rc:D+4:5}]}}
\end{lemma}

\begin{figure}[h!]
    \centering
    \includegraphics[width=0.25\columnwidth,page=7]{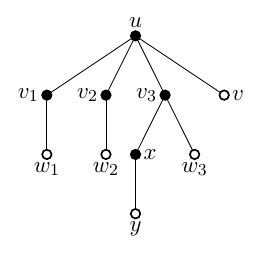}
    \caption{An illustration of Lemma~\ref{lem:D+4:4two2}.}
    \label{fig:lem:D+4:4two2}
\end{figure}

\begin{proof}
    Suppose to the contrary that $G$ has a $4$-vertex $u$ with two $2$-neighbors $v_1$, $v_2$ and one $3_1$-neighbor $v_3$.
We follow the labeling of the vertices as Figure~\ref{fig:lem:D+4:4two2}.
    Let $S=\{u, v_1, v_2, v_3, x\}$.
    Note that $x\neq v_i$ for each $i\in\{1,2\}$ by Lemma~\ref{lem:D+c}~(ii), and so $S$ has five distinct vertices.
    For each $i\in\{1,2\}$, note that $u\neq y$ and $v_3\neq w_i$ by Lemma~\ref{lem:D+c}~(ii).
    Also, $v_i\notin\{w_3,w_{3-i},y\}$, $x\notin\{w_1,w_2,v\}$ by Lemma~\ref{lem:D+c}~(iv).
    Thus $S\cap\{w_1,w_2,w_3,y,v\}=\emptyset$.
    Note that $|A_{\vp}(u)|\ge 1$, $|A_{\vp}(v_i)|\ge 3$ for each $i\in\{1,2\}$, $|A_{\vp}(v_3)|\ge 2$, and $|A_{\vp}(x)|\ge 3$ by \eqref{eq:Av}.

    First, we color $u$ with a color $\gamma$ from $A_{\vp}(u)$.
For simplicity, we define $A_i$ as follows.
Let $A_{i}:= (A_{\varphi }(v_{i})\cup \{\varphi (v)\}) \setminus \{\gamma \}$, if one of the following holds.
\begin{itemize}
\item $i\in\{1,2\}$, $w_i=w_j$ for some $j$ and $|A_{\varphi }(v_{i})|=4$,
\item $i=3$ and $|A_{\varphi }(v_{i})|=2$.
\end{itemize}
Let $A_{i}:=A_{\varphi }(v_{i}) \setminus \{\gamma \}$ otherwise.
For
each $i\in \{1,2,3\}$, 
since $\varphi (v)\notin A_{\varphi }(v_{i})$, it holds that  
$|A_{i}|\ge 2$. Let $A_{4}=\{\varphi (v)\}$. We color
$v_{1},v_{2},v_{3}$ as follows.

    If $w_1=w_2=w_3$, then $A_1=A_2$ and $A_3\subset A_1$, and so we color $v_i$'s with a same color $\beta\in A_3\setminus\{\vp(v)\}$.
    Suppose that it is not the case of $w_1=w_2=w_3$.
    Without loss of generality, let $w_2\neq w_3$.
    We take subsets $A'_2$ and $A'_3$ of $A_2$ and $A_3$, respectively, of size exactly two.
    If $w_1=w_j$ for some $j\in\{2,3\}$, then $|A_1| \ge 4$ and so we take a subset $A'_1$ of $A_1\setminus A'_j$ of size two.
    If $w_1 \neq w_j$ for each $j\in\{2,3\}$, then let $A'_1=A_1$.
    Now, we apply Lemma~\ref{lem:even} to obtain a system of odd representative $(\beta_1,\beta_2,\beta_3,\vp(v))$ for $(A'_1,A'_2,A'_3,A_4)$. We color $v_i$ with $\beta_i$ for each $i\in \{1,2,3\}$.
    Then the resulting coloring is an SO $(\D+4)$-coloring of $G-x$, which is a contradiction to Lemma~\ref{lem:G-x:2vtex}.
\end{proof}

We complete the proof of Theorem~\ref{thm:soD+4} by discharging technique.
Let $\mu(v):=\deg_G(v)$ be the initial charge of a vertex $v$. Then $\sum_{v\in V(G)}\mu(v) \le \frac{20|V(G)|}{7}$, since $mad(G)\le \frac{20}{7}$.
We let $\mu^*(v)$ be the final charge of $v$ after the following discharging rules:
\begin{enumerate}[(R1)]
    \item\label{4rule:1-1}
    Each $3^+$-vertex sends charge $\frac{3}{7}$ to each of its $2$-neighbors.
    \item\label{4rule:2-1}
    Each $3$-vertex without a $2$-neighbor or $4^+$-vertex sends charge $\frac{1}{7}$ to each of its $3_1$-neighbors.
\end{enumerate}

Take a vertex $v$.
By \textbf{[C\ref{rc:D+4:1}]}, $\deg_G(v)\ge 2$.
If $\deg_G(v)=2$, then $v$ has only $3^+$-neighbors by \textbf{[C\ref{rc:D+4:2}]}, and so
$\mu^*(v)=2+\frac{3}{7}\cdot 2=\frac{20}{7}$ by (R\ref{4rule:1-1}).
Suppose that $\deg_G(v)=3$. If $v$ has no $2$-neighbor, then $v$ has at most one $3_1$-neighbor by \textbf{[C\ref{rc:D+4:4}]}, so $\mu^*(v)\ge 3-\frac{1}{7}=\frac{20}{7}$ by (R\ref{4rule:2-1}).
Suppose that $v$ has a $2$-neighbor.
By \textbf{[C\ref{rc:D+4:2}]} and \textbf{[C\ref{rc:D+4:3}]}, $v$ has exactly one $2$-neighbor and two other $3^+$-neighbors{, which are either $4^+$-vertices or $3$-vertices having no $2$-neighbors}, and so $\mu^*(v)=3-\frac{3}{7}+\frac{1}{7}\cdot 2=\frac{20}{7}$ by (R\ref{4rule:1-1}) and (R\ref{4rule:2-1}).
Suppose that $\deg_G(v)=4$.
If $v$ has at most one $2$-neighbor, then $\mu^*(v)\ge 4-\frac{3}{7}-\frac{1}{7}\cdot 3>\frac{20}{7}$ by (R\ref{4rule:1-1}) and (R\ref{4rule:2-1}).
If $v$ has at least two $2$-neighbors, then
$v$ has exactly two $2$-neighbors and no $3_1$-neighbor by \textbf{[C\ref{rc:D+4:2}]} and \textbf{[C\ref{rc:D+4:5}]}, and so
 $\mu^*(v)= 4-\frac{3}{7}\cdot 2>\frac{20}{7}$
by (R\ref{4rule:1-1}).
If $\deg_G(v)\ge 5$, then it has at most $(\deg_G(v)-1)$ $2$-neighbors by \textbf{[C\ref{rc:D+4:add}]}, and so $\mu^*(v)\ge \deg_G(v)-\frac{3}{7}\cdot (\deg_G(v)-1) -\frac{1}{7}> \frac{20}{7}$ by (R\ref{4rule:1-1}) and (R\ref{4rule:2-1}).

From the fact that $\sum_{v\in V(G)}\mu(v)=\sum_{v\in V(G)}\mu^*(v) \le \frac{20|V(G)|}{7}$, we can conclude that every vertex has final charge exactly $\frac{20}{7}$. Then $\D=3$, and $V(G)=X_2 \cup X_{3_1}\cup Y$ if we let $X_2$ and $X_{3_1}$ be the set of $2$-, $3_1$-vertices of $G$, respectively, and let $Y$ be the set of $3$-vertices without $2$-neighbors.
Moreover, every $3$-vertex in $Y$ has exactly one $3_1$-neighbor, and so each connected component of $G[Y]$ is $2$-regular.
Lastly, note that both $X_2$ and $X_{3_1}$ are independent sets by \textbf{[C\ref{rc:D+4:2}]} and \textbf{[C\ref{rc:D+4:3}]}.

\begin{figure}[h!]
    \centering
    \includegraphics[width=0.33\textwidth,page=8]{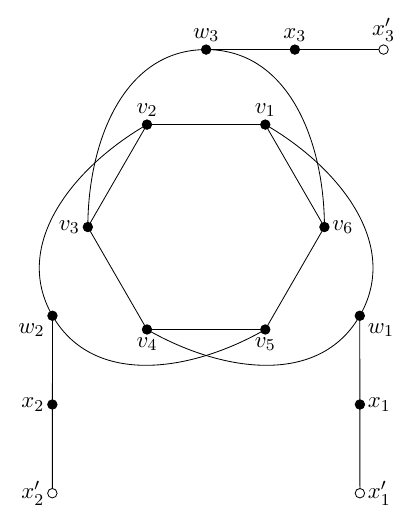}
    \caption{{An }illustration for the proof of Lemma~\ref{lem:Brooks}.}
    \label{fig:GwithX}
\end{figure}

Now, we finish the proof by showing that $G^2$ has a proper $7$-coloring, which is also an SO coloring of $G$.
Let $H:=G^2[Y]$, the subgraph of $G^2$ induced by $Y$.
Note that $\D(H)\le 5$, since each vertex of $H$ has at most four {neighbors} with distance at most two {in $G$} through only vertices in $Y$, and {has} one {neighbor} with distance at most two {in $G$} through a vertex not in $Y$.

\begin{lemma}\label{lem:Brooks}
There is no connected component in $H$ that is isomorphic to $K_6$.
\end{lemma}
\begin{proof}
Suppose to the contrary that there is a connected component $K$ in $H$ that is isomorphic to $K_6$.
Then the vertices of $K$ form a cycle in $G$, say $K:v_1v_2 v_3 v_4 v_5 v_6v_1$, and $v_i$ and $v_{i+3}$ have a common $3_1$-neighbor $w_i$ for each $i\in\{1,2,3\}$.
Note that $v_i$'s and $w_j$'s are distinct.
See Figure~\ref{fig:GwithX}.
For each $i\in \{1,2,3\}$, let $x_i$ be a 2-neighbor of $w_i$, $x'_i$ be the neighbor of $x_i$ other than $w_i$.
If $x_1=x_2$, then $x_3$ is a cut vertex of $G$, which is a contradiction by permuting the colors of an SO 7-coloring of a connected component of $G-x_3$.
Thus $x_1$, $x_2$, and $x_3$ are distinct. 

Let $S=V(K)\cup \{w_1,w_2,w_3,x_1,x_2,x_3\}$.
Note that $|A_{\vp}(w_i)|\ge 6$, $|A_{\vp}(x_i)|\ge 4$ for each $i\in\{1,2,3\}$, $|A_{\vp}(v_j)|\ge 7$ for each $j\in\{1,\ldots, 6\}$ by \eqref{eq:Av}.
First, we color $w_i$ with a color in $A_{\vp}(w_i)$ for each $i\in\{1,2,3\}$, $v_i$ with a color in $A_{\vp}(v_j)$ for each $j\in\{1,3,4,5\}$ so that
 the colors of seven vertices $w_i$'s and $v_j$'s are all distinct.
Then we color $v_2$ and $v_6$ with the same color with $w_1$.
We denote this resulting coloring of $G-\{x_1,x_2,x_3\}$ as $\vp$ again.
For each $i\in\{1,2,3\}$, $|A_{\vp}(x_i)|\ge 1$ by \eqref{eq:Av} and
we color $x_i$ with a color in $A_{\vp}(x_i)$.
It gives an SO $7$-coloring of $G$.
\end{proof}

From Lemma~\ref{lem:Brooks} and Theorem~\ref{thm:brooks}, there is a proper $5$-coloring $\phi_H:V(H)\rightarrow \{1,\ldots,5\}$.
We color each vertex in $X_{3_1}$ with $6$ or $7$ so that two vertices having a common $2$-neighbor have distinct colors, and then color each vertex $x\in X_2$ with a color not appeared in $N_{G^2}(x)$.
{It} is possible since $|N_{G^2}(x)|\le 6$. Then it gives a proper coloring of $G^2$ and completes the proof.

\section{Strong odd $(\D+3)$-coloring}

Like as the previous section, in each proof, we always start with defining a nonempty set $S\subset V(G)$ and $\vp$ is always an SO $(\Delta(G)+{3})$-coloring $\vp$ of $G'=G-S$. 

Let $G$ be a minimal counterexample to Theorem~\ref{thm:so6coloring}. 
Let $\Delta:=\Delta(G)$. 
If $\D\le 2$, then it is easy to see that $G$ has an $(\Delta(G)+{3})$-coloring. 
Thus $\D\ge 3$. We will show that the following do not appear in $G$.
For simplicity, let  \[Z=\{v\in V(G)\mid v\text{ is a }3_1\text{-vertex with two $3_1$-neighbors}\}.\]
\begin{enumerate}[\bfseries {[}C1{]}]
    \item\label{rc:D+3:1-1} A $1^-$-vertex. \vspace{-0.25cm}
    \item\label{rc:D+3:2-1} A $d_{d-1}$-vertex for all $2\le d\le 3$. \vspace{-0.25cm}
    \item\label{rc:D+3:3-1} A $d_d$-vertex for all $d\ge 4$. \vspace{-0.25cm}
    \item\label{rc:D+3:5:NEW} A $2$-vertex with only $3$-neighbors in which at least one  is in $Z$.  \vspace{-0.25cm}
    \item\label{rc:D+3:6:NEW} A $4_3$-vertex with a {$3$}-neighbor or a $4_3$-neighbor.
\end{enumerate}
Note that \textbf{[C\ref{rc:D+3:1-1}]}-\textbf{[C\ref{rc:D+3:3-1}]} do not appear in $G$ by Lemma~\ref{lem:D+c}~(i) and (iv).
We will prove \textbf{[C\ref{rc:D+3:5:NEW}]} does not exist in $G$ by considering two cases depending on $\D$.

\begin{lemma}\label{lem:4.1}
When $\Delta \ge 4$, the graph $G$ has no  $2$-vertex  with two $3$-neighbors in which one  is in $Z$. {\rm\textbf{[C\ref{rc:D+3:5:NEW}]}}
\end{lemma}
\begin{figure}[h!]
    \centering
    \includegraphics[width=0.3\columnwidth,page=11]{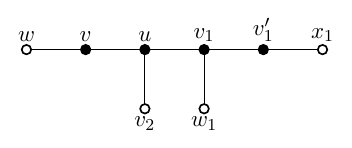}
     \caption{An illustration of Lemma~\ref{lem:4.1}.}
    \label{fig:lem:4.1}
\end{figure}
\begin{proof} Suppose to the contrary that there is a $2$-vertex $v$ with $3$-neighbors $u,w$ and $u\in Z$. Let $v_1,v_2$ be the $3_1$-neighbors of $u$, $v'_1$ be the 2-neighbor of $v_1$, $w_1$ be the other neighbor of $v_1$, and $x_1$ be the neighbor of $v_1'$ other than $v_1$. See Figure~\ref{fig:lem:4.1}.
If we have an SO $(\D+3)$-coloring of $G-v$ such that the
colors of $u$ and $w$ are distinct, then it is a contradiction to Lemma~\ref{lem:G-x:2vtex}, since $v$ always has an available color from the facts that $|\mathcal{C}|=\D+3\ge 7$.

Let $S=\{u,v,v_1,v'_1\}$.
By Lemma~\ref{lem:D+c}~(ii), $v\neq v'_1$, and so $S$ has four distinct vertices.
By Lemma~\ref{lem:D+c}~(iv), $u\neq x_1$, $v\neq w_1$, $v\neq x_1$, $v'_1\neq w$ and so $S\cap \{w,v_2,w_1,x_1\}=\emptyset$. 
Also $w\neq v_2$ and $w_1\neq x_1$ by Lemma~\ref{lem:D+c}~(ii).
Recall that $v_2$ is a $3$-vertex and so note that 
$|A_{\vp}(u)|\ge\D-2\ge 2$, $|A_{\vp}(v_1)|\ge 1$, and  $|A_{\vp}(v'_1)|\ge 2$ by \eqref{eq:Av}.
If we can color $x$ with a color in $A_{\vp}(x)$ for each $x\in\{u,v_1,v'_1\}$ so that their colors are distinct, then it is an SO $(\D+3)$-coloring of $G-v$ such that the colors of $u$ and $w$ are distinct, which contradicts to Lemma~\ref{lem:G-x:2vtex}.
Thus we may assume that $\D=4$, $A_{\vp}(u)=A_{\vp}(v'_1)=\{1,2\}$,  $A_{\vp}(v_1)\subset\{1,2\}$ and $\vp(w_1)=3$. 
Then $3\not\in \vp(N_G[x_1])\cup \vp(N_G[v_2])$, $\vp(w)\neq 3$.
We color $v'_1$ and $u$ with the color $3$ and $v_1$ with a color in $A_{\vp}(v_1)$. 
It is an SO $7$-coloring of $G-v$ such that the colors of $u$ and $w$ are distinct, which contradicts to Lemma~\ref{lem:G-x:2vtex}.
\end{proof}

\begin{lemma}\label{lem:D+3:no3cycle}
When $\Delta=3$, the graph $G$ has no triangle consisting of $3_1$-vertices.
\end{lemma}

\begin{figure}[h!]
    \centering
    \includegraphics[width=0.23\columnwidth,page=9]{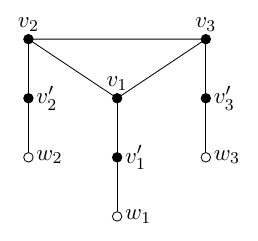}
    \caption{An illustration of Lemma~\ref{lem:D+3:no3cycle}.}
    \label{fig:lem:D+3:no3cycle}
\end{figure}

\begin{proof}
Suppose to the contrary that $G$ has a triangle $v_1v_2v_3v_1$ consisting of only $3_1$-vertices.
Let $v'_i$ be the $2$-neighbor of $v_i$ and $N_G(v'_i)=\{v_i, w_i\}$ for each $i\in\{1,2,3\}$.
See Figure~\ref{fig:lem:D+3:no3cycle}.
Let $S=\{v_1,v_2,v_3,v'_1,v'_2,v'_3\}$. 
By Lemma~\ref{lem:D+c}~(ii), $v'_1$, $v'_2$, $v'_3$ are distinct, and so $S$ has six distinct vertices.

By \eqref{eq:Av},  $|A_{\vp}(v_i)|\ge \Delta+2=5$, and $|A_{\vp}(v'_i)|\ge 3$ for each $i\in\{1,2,3\}$. 
If we can choose  a color in $A_{\vp}(x)$ for each $x\in S$ to color with distinct colors, then it results in an SO $(\Delta+3)$-coloring of $G$. Thus $A_{\vp}(v_1)=A_{\vp}(v_2)=A_{\vp}(v_3)$ and so $\vp(w_1)=\vp(w_2)=\vp(w_3)$.  
 Since $|A_{\vp}(v'_i)|\ge 3$ and $\vp(w_1)=\vp(w_2)$, it follows that $\gamma\in A_{\vp}(v'_1)\cap  A_{\vp}(v'_2)\neq\emptyset$ for some $\gamma$.
We color $v'_1$ and $v'_2$ with $\gamma$.
For $x\in N_G[v_3]$,
we color $x$ with a color in $A_\vp(x)\setminus\{\gamma\}$ so that all colors of the vertices in $N_G[v_3]$ are distinct.
It is an SO $6$-coloring of $G$.
\end{proof}

\begin{lemma}\label{lem:two3_1}
When $\Delta=3$, $Z=\emptyset$. That is, if $G$ is subcubic, then 
$G$ has no $3_1$-vertex with two $3_1$-neighbors
{\rm\textbf{[C\ref{rc:D+3:5:NEW}]}}
\end{lemma}
\begin{figure}[h!]
    \centering
    \includegraphics[width=0.25\columnwidth,page=10]{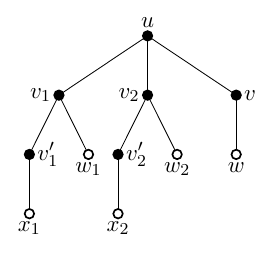}
     \caption{An illustration of Lemma~\ref{lem:two3_1}.}
    \label{fig:lem:two3_1}
\end{figure}
\begin{proof}
Suppose to the contrary that $G$ has a $3_1$-vertex $u$ with two $3_1$-neighbors. 
We follow the labeling of the vertices as Figure~\ref{fig:lem:two3_1}.
Note that $v'_1$, $v'_2$, $w$, $w_1$, and $w_2$ are distinct since $G$ has no $4$-cycle.
Let $S=\{u,v,v_1,v_2,v'_1,v'_2\}$.
By Lemma~\ref{lem:D+c}~(ii), $v\neq v'_i$ for each $i\in\{1,2\}$, and so $S$ has six distinct vertices.
Note that for each $i\in\{1,2\}$, $u\neq x_i$, $v_i\neq w$ by Lemma~\ref{lem:D+c}~(ii).
Also, $v\notin\{x_1,x_2,w_1,w_2\}$, $v'_i\notin\{x_1,x_2,w_1,w_2,w\}$, $v_i\neq x_{3-i}$ by Lemma~\ref{lem:D+c}~(iv), and $v_i\neq w_{3-i}$ by Lemma~\ref{lem:D+3:no3cycle}.
Thus $S\cap \{x_1, x_2, w_1, w_2, w\}= \emptyset$. 

In addition, $w$, $x_1$, $x_2$ are distinct by Lemma~\ref{lem:D+c}~(iv).
Note that $|A_{\vp}(v_i)|\ge 2$ for each $i\in\{1,2\}$ by \eqref{eq:Av}.
Also, since $w$ is a $2$-vertex in $G'$, $|A_{\vp}(v)|\ge 3$ by \eqref{eq:Av}.
We now proceed to prove $5$ claims.

\begin{claim}\label{clm:nocommon}
$A_{\vp}(v)\cap A_{\vp}(v_1) \cap A_{\vp}(v_2)= \emptyset$.
\end{claim}
\begin{proof}
Suppose to the contrary that $\alpha \in A_{\vp}(v)\cap A_{\vp}(v_1) \cap A_{\vp}(v_2)$.
We give a color $\alpha$ to $v$, $v_1$, and $v_2$, and denote this coloring of $G-\{u,v'_1,v'_2\}$ by $\vp$ again.
Note that $|A_{\vp}(v'_i)|\ge 1$ for each $i\in\{1,2\}$ by \eqref{eq:Av}.
Next, we color $v'_i$ with a color $\beta_i$ in $A_{\vp}(v'_i)$ for each $i\in\{1,2\}$.

If there is a color in $\mathcal{C} \setminus \{\vp(w_1),\vp(w_2),\vp(w),\alpha,\beta_1,\beta_2\}$, then we color $u$ with the color and it is done.
Thus $\{\vp(w_1),\vp(w_2),\vp(w),\alpha,\beta_1,\beta_2\}=\mathcal{C}$.
Moreover, $A_{\vp}(v'_i)=\{\beta_i\}$ and so
$\vp(w_i)\notin \vp(N_{G'}[x_i])$ for each $i\in\{1,2\}$.
Coloring both $v'_1$ and $u$ with $\vp(w_1)$ gives an SO $6$-coloring of $G$.
This concludes the proof of the claim.
\end{proof}

\begin{claim}\label{clm:colorofv}
There is a color $\alpha_i\in A_{\vp}(v_i)$ for each $i\in\{1,2\}$ such that $|A_{\vp}(v)\setminus\{\alpha_1,\alpha_2\}|\ge 2$ and $\alpha_1\neq \alpha_2$.
\end{claim}
\begin{proof}
Recall that $|A_{\vp}(v_1)|,|A_{\vp}(v_2)|\ge 2$ and $|A_{\vp}(v)|=3$.  
Suppose that $A_{\vp}(v_i) \subset A_{\vp}(v)$ for each $i\in\{1,2\}$. Then  $|A_{\vp}(v_1)\cap A_{\vp}(v_2)|\ge |A_{\vp}(v_1)|+|A_{\vp}(v_2)|-|A_{\vp}(v)|\ge 2+2-3=1$. An element in $A_{\vp}(v_1)\cap A_{\vp}(v_2)$ belongs to $A_{\vp}(v_1)\cap A_{\vp}(v_2)\cap A_{\vp}(v) $,  a contradiction to Claim~\ref{clm:nocommon}.
Thus  $A_{\vp}(v_i) \setminus A_{\vp}(v)\neq \emptyset$ for some $i\in\{1,2\}$. We may assume that $i=1$. We take $\alpha_1\in A_{\vp}(v_1) \setminus A_{\vp}(v)$ and $\alpha_2\in A_{\vp}(v_2)\setminus\{\alpha_1\}$, and so $|A_{\vp}(v)\setminus\{\alpha_1,\alpha_2\}|\ge 2$.
This concludes the proof of the claim.
\end{proof}

We color $v_i$ with $\alpha_i$  for each $i\in\{1,2\}$, where $\alpha_1$ and $\alpha_2$ satisfy Claim~\ref{clm:colorofv}, and denote this coloring of $G-\{u,v,v'_1,v'_2\}$ by $\vp$ again. Then for each $i\in\{1,2\}$,
$|A_{\vp}(v'_i)|\ge 1$  by \eqref{eq:Av}, and
\[(\S)\qquad \text{if }|A_{\vp}(v'_i)|=1\text{ then }\vp(w_i) \notin \vp(N_{G'} [x_i]). \qquad \]

\begin{claim}\label{claim:last:v}
There is no  SO $6$-coloring of $G-v$ that is
an extension of $\vp$ such that the color of $u$ is not $\vp(w)$.
\end{claim}
\begin{proof} Suppose that we color $v'_1$, $v'_2$, and $u$ to obtain an SO $6$-coloring $\vp^*$ of $G-v$ satisfying $\vp^*(u)\neq \vp(w)$.
By Claim~\ref{clm:colorofv}, there is at least one available color in $A_{\vp^*}(v)$. It contradicts to Lemma~\ref{lem:G-x:2vtex}.
This concludes the proof of the claim.
\end{proof}

Let $B=A_{\vp}(v'_1) \times A_{\vp}(v'_2)$.
Note that $B\neq \emptyset$.
If there is a color $\gamma\in A_{\vp}(u)\setminus\{ \beta_1,\beta_2\}$
for some $(\beta_1,\beta_2)\in B$, then coloring $v'_1$, $v'_2$ and $u$ with $\beta_1$, $\beta_2$ and $\gamma$, respectively, gives an SO $6$-coloring of $G-v$ such that $\gamma\neq \vp(w)$, which contradicts to Claim~\ref{claim:last:v}.
Thus $A_{\vp}(u)\setminus\{ \beta_1,\beta_2\}=\emptyset$, and so
\begin{eqnarray}\label{eq:C}
\mathcal{C}&=&\{\alpha_1,\alpha_2,\beta_1,\beta_2,\vp(w),\vp(w_1),\vp(w_2)\}  \qquad\text{ for every } (\beta_1,\beta_2)\in B.
\end{eqnarray}
Note that  for every $(\beta_1,\beta_2)\in B$, $\beta_i\neq \vp(w_i)$ by the definition of $A_{\vp}(v'_i)$ for each $i\in\{1,2\}$ and $\{\beta_i, \vp(w_i)\} \cap \{\alpha_1,\alpha_2,\vp(w)\}=\emptyset$ for some $i\in\{1,2\}$.

\begin{claim}\label{clm:base}
For every $(\beta_1,\beta_2)\in B$, it holds that $|\{\beta_1,\beta_2,\vp(w_1),\vp(w_2)\}|=3$ and $A_{\vp}(v'_1)\cup A_{\vp}(v'_2)\subset\{\beta_1,\beta_2,\vp(w_1),\vp(w_2)\}$.
\end{claim}
\begin{proof} By \eqref{eq:C},    $|\{\beta_1,\beta_2,\vp(w_1),\vp(w_2)\}|\ge 3$ for every $(\beta_1,\beta_2)\in B$.
Suppose to the contrary $|\{\beta_1,\beta_2,\vp(w_1),\vp(w_2)\}|=4$ for some $(\beta_1,\beta_2)\in B$.
Without loss of generality, let $\{\beta_1, \vp(w_1)\} \cap \{\alpha_1,\alpha_2,\vp(w)\}=\emptyset$.
If $|A_{\vp}(v'_1)|\ge 2$, then we color $v'_1$, $v'_2$, and $u$ with a color with $\beta'_1\in A_{\vp}(v'_1)\setminus\{\beta_1\}$, $\beta_2$, and $\beta_1$, respectively, which gives an SO $6$-coloring of $G-v$ such that the color of $u$ is not $\vp(w)$, a contradiction to Claim~\ref{claim:last:v}.
Thus $|A_{\vp}(v'_1)|=1$  and so by ($\S$), $\vp(w_1)\notin \vp(N_{G'}[x_1])$.
Then we color $v'_1$, $v'_2$, and $u$ with $\vp(w_1)$, $\beta_2$, and $\vp(w_1)$, respectively, which gives an SO $6$-coloring of $G-v$.
Then the color of $u$ is not $\vp(w)$ and so it is a contradiction to Claim~\ref{claim:last:v}.
Thus $|\{\beta_1,\beta_2,\vp(w_1),\vp(w_2)\}|=3$ and so
$\{\alpha_1,\alpha_2,\vp(w)\}\cap \{\beta_1,\beta_2,\vp(w_1),\vp(w_2)\}=\emptyset$ for every $(\beta_1,\beta_2)\in B$. Therefore each $A_{\vp}(v'_i)$ is a subset of $ \{\beta_1,\beta_2,\vp(w_1),\vp(w_2)\}$.
This concludes the proof of the claim.
\end{proof}

If $\vp(w_1)=\vp(w_2)$, then for each $i\in \{1,2\}$, $|A_{\vp}(v'_i)|=1$ by Claim~\ref{clm:base} and so $\vp(w_1)\notin \vp(N_{G'}[x_{i}])$ by ($\S$).
Then we color $v'_1$, $v'_2$, and $u$ with $\vp(w_1)$, and it is an SO $6$-coloring of $G-v$ such that the color of $u$ is not $\vp(w)$, which contradicts to Claim~\ref{claim:last:v}. 
Thus $\vp(w_1)\neq \vp(w_2)$.
Together with \eqref{eq:C} and Claim~\ref{clm:base}, \[ |\{\alpha_1,\alpha_2,\vp(w),\vp(w_1),\vp(w_2)\}|=5.\] 
Let $\beta^*$ be the color not in $\mathcal{C}\setminus\{ \alpha_1,\alpha_2,\vp(w),\vp(w_1),\vp(w_2) \}$.
In addition, for each $i\in\{1,2\}$,
$|A_\vp(v'_i)|\le 2$, and moreover, 
$A_\vp(v'_i)\subset\{\beta^*,\vp(w_{3-i})\}$.

By Claim~\ref{clm:base}, since $(\vp(w_2),\vp(w_1) )\notin B$,
\[ (\ddag) \qquad \vp(w_2)\not\in A_{\vp}(v'_1)\text{ or }\vp(w_1)\not\in A_{\vp}(v'_2).\]
\begin{claim}\label{clm:d1d2}
For some $i\in\{1,2\}$,  $\vp(w_i)\not\in \vp(N_{G'}[x_i])$ and  $\beta^*\in A_\vp(v'_{3-i})$.
\end{claim}
\begin{proof} 
Note that $\beta^*\notin A_\vp(v'_{3-i})$ implies that 
$A_{\vp}(v'_{3-i})=\{\vp(w_i)\}$, since $A_\vp(v'_{3-i})\neq \emptyset$. 
Suppose to the contrary that each $\vp(w_i)$ satisfies (P) or (Q), where (P) and (Q) are properties defined as follows:
{\rm{(P)}} $\vp(w_i)\in \vp(N_{G'}[x_i])$ and  {\rm{(Q)}}  {$A_{\vp}(v'_{3-i})=\{\vp(w_i)\}$}.

By ($\ddag$), $\vp(w_i)$ does not satisfy (Q) for some $i\in\{1,2\}$. Then $\vp(w_i)$ satisfies (P), and so $|A_{\vp}(v'_i)|\ge 2$ by ($\S$). Then $\vp(w_{3-i})$ does not satisfy (Q). By the assumption, $\vp(w_{3-i})$ also satisfies (P) and so  $|A_{\vp}(v'_{3-i})|\ge2$ by ($\S$). 
Therefore, $|A_{\vp}(v'_{1})|=|A_{\vp}(v'_{2})|=2$, and thus
  $A_\vp(v'_1)=\{\beta^*,\vp(w_{2})\}$ and $A_\vp(v'_{2})=\{\beta^*,\vp(w_1)\}$, which is a contradiction to ($\ddag$).
\end{proof}

By Claim~\ref{clm:d1d2}, without loss of generality, we may let $\vp(w_1)\notin \vp(N_{G'}[x_1])$ and $\beta^*\in A_\vp(v'_2)$.
We color $v'_1$, $v'_2$, and $u$ with $\vp(w_1)$, $\beta^*$, and $\vp(w_1)$, respectively, since $\vp(w_1)\neq \vp(w_2)$.
It gives an SO $6$-coloring of $G-v$ and the color of $u$ is not $\vp(w)$, which is a contradiction to Claim~\ref{claim:last:v}.
\end{proof}

\begin{lemma}\label{lem:4.2}
The graph $G$ has no $4_3$-vertex with a {$3$}-neighbor or a $4_3$-neighbor. {\rm\textbf{[C\ref{rc:D+3:6:NEW}]}}
\end{lemma}
\begin{figure}[h!]
    \centering
    \includegraphics[width=0.35\columnwidth,page=12]{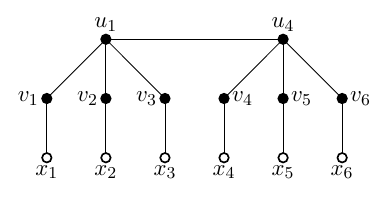}
     \caption{An illustration of Lemma~\ref{lem:4.2}.}
    \label{fig:lem:4.2}
\end{figure}
\begin{proof}
A $4_3$-vertex does not have a {$3$}-neighbor by Lemma~\ref{lem:D+c}~(iv). Also, $\D=4$ by Lemma~\ref{lem:D+c}~(iv) again.
Let $u_1$ and $u_4$ be adjacent $4_3$-vertices. Let $N_G(u_i)=\{u_{5-i}, v_i,v_{i+1},v_{i+2}\}$ for each $i\in\{1,4\}$, and 
$x_i$ be the neighbor of $v_i$ other than $u_1,u_4$ for each $i\in\{1,\ldots,6\}$. See Figure~\ref{fig:lem:4.2}. By Lemma~\ref{lem:D+c}~(ii), all $v_i$'s are distinct. 
 {Also $v_i\neq x_j$ for $i,j\in \{1,\ldots, 6\}$ by Lemma~\ref{lem:D+c}~(iv), and $x_i,x_{i+1},x_{i+2}$ are all distinct for $i\in\{1,4\}$ since $G$ is $C_4$-free.}
We begin with the following claim.

\begin{claim}\label{claim:4_3} For each $i\in\{1,4\}$,
 there is no SO $7$-coloring $\vp$ of $G-\{v_i,v_{i+1},v_{i+2}\}$ such that the color of $u_i$ is not equal to colors of $x_i,x_{i+1},x_{i+2}$.
\end{claim}
\begin{proof}
Suppose that there is an SO $7$-coloring of $G-\{v_1,v_{2},v_3\}$ satisfying the condition of the claim. 
 {Note that $|A_{\vp}(v_i)|\ge 1$ for each $i\in \{1,2,3\}$ by \eqref{eq:Av}.}
Let {$A_i:=A_{\vp}(v_i)\cup \{\vp(u_4)\}$ if $|A_{\vp}(v_i)|=1$, and $A_i:=A_{\vp}(v_i)$ otherwise} for each $i\in \{1,2,3\}$, and let $A_4=\{ \vp(u_4) \}$. 
Note that $|A_i|\ge 2$ for each $i\in \{1,2,3\}$. 
Applying Lemma~\ref{lem:even}, there is a system of odd representative ($\beta_1,\beta_2,\beta_3,\vp(u_4)$) for $(A_1,A_2,A_3,A_4)$. We color each $v_i$ with $\beta_i$, which is an SO $7$-coloring of $G$.
 {The argument is the same with $G-\{v_4,v_5,v_6\}$.}
This concludes the proof of the claim.
\end{proof}

Let $S=N_G(u_1)\cup N_G(u_4)$.  Note that $S$ has eight distinct vertices.
Also, $|A_{\vp}(v_i)|\ge 3$ and
$|A_{\vp}(u_j)|\ge  4$ for each $i\in\{1,2,3,4,5,6\}$ and $j\in\{1,4\}$ by \eqref{eq:Av}.
By Claim~\ref{claim:4_3}, we cannot color $x$ with a color in $A_{\vp}(x)$ for each $x\in N_G[u_4]$ so that the colors of $N_G[u_4]$ are distinct. 
Thus we may assume that $A_{\vp}(u_1)=A_{\vp}(u_4)=\{1,2,3,4\}$ and $A_{\vp}(v_j)\subset \{1,2,3,4\}$ for each $j\in\{4,5,6\}$.
Then $\gamma \in A_{\vp}(v_4)\cap A_{\vp}(v_5)\cap A_{\vp}(v_6)$ for some $\gamma$ and so we color $v_4$, $v_5$, $v_6$ with $\gamma$, and then color $u_1$ and $u_4$ with distinct colors in $\{1,2,3,4\}\setminus\{\gamma\}$. It is an SO $7$-coloring of $G-\{v_1,v_{2},v_3\}$ {such that the color of $u_1$ is not equal to colors of $x_1, x_2, x_3$}, which contradicts to Claim~\ref{claim:4_3}.
\end{proof}

We complete the proof of Theorem~\ref{thm:so6coloring} by discharging technique.
Let $\mu(v):=\deg_G(v)$ be the initial charge of a vertex $v$, and let $\mu^*(v)$ be the final charge of $v$ after the  discharging rules:

\begin{enumerate}[(R1)]
    \item \label{3Xrule:1}
    Each $3$-vertex in $Z$ sends charge $\frac{3}{11}$ to each of its $2$-neighbors.
    \item \label{3notXrule:2}
    Each $3$-vertex not in $Z$ sends charge $\frac{4}{11}$ to each of its $2$-neighbors.
    \item \label{4rule:1}
    Each $4^+$-vertex sends charge $\frac{5}{11}$ to each of its $2$-neighbors.
    \item \label{4_3rule:2}
    Each $4$-vertex with at most two $2$-neighbors sends charge $\frac{1}{11}$ to each of its $3_1$-neighbors and $4_3$-neighbors.
    \item \label{5rule:1}
    Each $5^+$- or $3$-vertex without $2$-neighbor sends charge $\frac{1}{11}$ to each of its $3_1$-neighbors and $4_3$-neighbors.
\end{enumerate}

Take a vertex $v$. We will show that $\mu^*(v) \ge \frac{30}{11}$.
By \textbf{[C\ref{rc:D+3:1-1}]}, $\deg_G(v)\ge 2$.
Suppose that  $\deg_G(v)=2$. Then $v$ has only $3^+$-neighbors by \textbf{[C\ref{rc:D+3:2-1}]}.
If $v$ has no neighbor in $Z$,
then $\mu^*(v)\ge 2+\frac{4}{11}\cdot 2=\frac{30}{11}$ by (R\ref{3notXrule:2}) and (R\ref{4rule:1}).
If $v$ has a neighbor in $Z$, then by \textbf{[C\ref{rc:D+3:5:NEW}]} $v$ also has a $4^+$-neighbor and so $\mu^*(v)\ge 2+\frac{3}{11}+\frac{5}{11}=\frac{30}{11}$ by (R\ref{3Xrule:1}) and (R\ref{4rule:1}).

Suppose $\deg_G(v)=3$.
By \textbf{[C\ref{rc:D+3:2-1}]}, it has at most one $2$-neighbor.
If $v$ has no $2$-neighbor, then $\mu^*(v)\ge 3-\frac{1}{11}\cdot 3=\frac{30}{11}$ by (R\ref{5rule:1}).
Suppose $v$ is a $3_1$-vertex.
If $v\in Z$, then  by (R\ref{3Xrule:1}),  $\mu^*(v)=3-\frac{3}{11}=\frac{30}{11}$.
Suppose that $v\not\in Z$. Then $v$ has at least one $3^+$-neighbor $u$ that is not a $3_1$-vertex.
Then $u$ is either a $3$-vertex without $2$-neighbors, a $4$-vertex with at most two 2-neighbors by \textbf{[C\ref{rc:D+3:6:NEW}]}, or a $5^+$-vertex.
By (R\ref{4_3rule:2}) and (R\ref{5rule:1}), $\mu^*(v)\ge 3-\frac{4}{11}+\frac{1}{11}=\frac{30}{11}$.

Suppose that $\deg_G(v)=4$.
Then it has at most three $2$-neighbors by \textbf{[C\ref{rc:D+3:3-1}]}.
If $v$ has at most two $2$-neighbors, then
by (R\ref{4rule:1}) and (R\ref{4_3rule:2}), $\mu^*(v)\ge 4-\frac{5}{11}\cdot 2-\frac{1}{11}\cdot 2 >\frac{30}{11}$.
Suppose that $v$ has three $2$-neighbors. Let $u$ be the $3^+$-neighbor of $v$. Then $u$ is neither a {$3$}-vertex nor a $4_3$-vertex by \textbf{[C\ref{rc:D+3:6:NEW}]} and so $u$ sends $\frac{1}{11}$ to $v$ by {(R\ref{4_3rule:2}) and} (R\ref{5rule:1}).
Then $\mu^*(v)\ge 4-\frac{5}{11}\cdot 3+\frac{1}{11} =\frac{30}{11}$.
If $\deg_G(v)\ge 5$, then it has at most $(\deg_G(v)-1)$ $2$-neighbors by {\rm\textbf{[C\ref{rc:D+3:3-1}]}}, and so $\mu^*(v)\ge \deg_G(v)-\frac{5}{11}\cdot (\deg_G(v)-1)-\frac{1}{11}> \frac{30}{11}$ by (R\ref{4rule:1}) and (R\ref{5rule:1}). 
Thus $\mu^*(v) \ge \frac{30}{11}$ for each vertex $v$.
{Since $\sum \mu^*(v)=\sum \mu(v)\le \frac{30|V(G)|}{11}$, we have $\mu^*(v)=\frac{30}{11}$ for each vertex $v$.
Thus $v$ is either $2$-, $3$-, or $4_3$-vertex.
If $\D\ge 4$, then there is a $4_3$-vertex $x$ in $G$.
By {\rm\textbf{[C\ref{rc:D+3:3-1}]}}, there is a $3^+$-neighbor $y$ of $x$. 
However, $y$ is neither a $4_3$-vertex nor a $3$-vertex by {\rm\textbf{[C\ref{rc:D+3:6:NEW}]}}, which is a contradiction.}
Thus $\D=3$.
Moreover by Lemma~\ref{lem:two3_1}, note that $Z=\emptyset$.

Let $X_0$ be the set of $2$-vertices,
$X_1$ be the the set of $3_1$-vertices, and
$X_2$ be  the the set of $3$-vertices not in $X_1$.
The final charge of each vertex is exactly $\frac{30}{11}$, and it implies that $G[X_1]$ forms an induced matching, and $X_2$ forms an independent set. We color each $2$-vertex with a color 1, each $3$-vertex in $X_2$ with a color $2$. We will color the vertices of $X_1$ with colors $3,4,5,6$ as follows.
Let $H=G^2[X_1]$, subgraph of $G^2$ induced by $X_1$. Then $\D(H)\le 4$.
Moreover, by the structure of $G$, 
$H$ cannot have $K_5$ as a connected component, since  if a vertex $v$ in $X_1$ has four neighbors $v_1$, $v_2$, $v_3$, $v_4$ in $H$ and we let  $vv_1\in E(G)$ and $v$, $v_2$ have a common neighbor $w\in X_0$ in $G$, then in the graph $H$, $v_2$ cannot be adjacent to both $v_3$ and $v_4$.
By Theorem~\ref{thm:brooks}, there is a proper coloring $\vp_H$ of $H$ with the color set $\{3,4,5,6\}$, and we color each vertex $x\in X_1$ with the color $\vp_H(x)$.
Then one can check that it is an SO coloring of $G$, a contradiction.

\section{Remarks}

Regarding to Theorem~\ref{thm:soD+4}, we know that the bound of the strong chromatic number of subcubic planar graph is tight by a graph in Figure~\ref{fig:planarcubic7}. 
But we believe the bound in Theorem~\ref{thm:soD+4} is not tight if a graph is large enough. 

In odd coloring, Petru\v sevski and \v Skrekovski~\cite{petr2022odd} asked if  $\chi_o(G)\le 5$ for a planar graph $G$.
It is proved that if $G$ is a planar graph, then $\chi_o(G)\le 8$ by Petr and Portier~\cite{petr2022odd}.
We can also ask it is possible to limit the strong odd chromatic number of all planar graphs to a constant.
After we uploaded this article online as a preprint, \cite{caro2025strongodd} answered this question affirmatively.

\section*{Acknowledgements}
The authors thank the referee for helpful comments and suggestions that improved the paper.
Boram Park was supported by the National Research Foundation of Korea(NRF) grant funded by the Korea government(MSIT) (No. RS-2025-00523206), and supported by the New Faculty Startup Fund from Seoul National University. 
Hyemin Kwon supported by a KIAS Individual Grant (CG099301) at Korea Institute for Advanced Study.


\end{document}